\documentclass[a4paper,11pt,fleqn]{article}
\usepackage[final]{showkeys}
\usepackage[obeyspaces,hyphens,spaces]{url}
\renewcommand*\showkeyslabelformat[1]{%
\fbox{\parbox[t]{1.2 cm}{\raggedright\normalfont\tiny%
\url{#1}}}}

\usepackage[dvipsnames]{xcolor}
\usepackage{enumitem}
\usepackage{amsmath}
\usepackage{amssymb}
\usepackage{amsfonts}
\usepackage{theorem}
\usepackage{mathtools}
\usepackage{euscript}
\usepackage{exscale,relsize}
\usepackage{epic,eepic}
\usepackage{makeidx}
\usepackage{srcltx} 
\usepackage{charter}
\usepackage{etoolbox}
\usepackage{mathrsfs} 
\usepackage{authblk}
\newcommand{\email}[1]{\href{mailto:#1}{\nolinkurl{#1}}}
\usepackage[normalem]{ulem} 
\allowdisplaybreaks
\definecolor{labelkey}{rgb}{0,0.08,0.45}
\definecolor{refkey}{rgb}{0,0.6,0.0}
\definecolor{dblue}{rgb}{0,0.08,0.75}
\definecolor{Dblue}{HTML}{8602DC}
\definecolor{dgreen}{HTML}{025940}
\definecolor{dred}{HTML}{D90404}
\usepackage{upref}
\usepackage{hyperref}
\hypersetup{colorlinks=true,linktocpage=true,linkcolor=dblue,%
citecolor=dgreen,urlcolor=dred}
\makeatletter
\let\@msm@th@eqref\eqref
\renewcommand{\eqref}[1]{%
\begingroup
\leavevmode
\color{dblue}%
\@msm@th@eqref{#1}%
\endgroup
}
\makeatother

\renewcommand{\leq}{\ensuremath{\leqslant}}
\renewcommand{\geq}{\ensuremath{\geqslant}}

\newcommand{\scal}[2]{{\langle{{#1}\mid{#2}}\rangle}}

\newcommand{\menge}[2]{\big\{{#1}~\big|~{#2}\big\}} 
\newcommand{\Menge}[2]{\Big\{{#1}~\big|~{#2}\Big\}} 

\newcommand{\HH}{\ensuremath{{\mathcal H}}}

\newcommand{\emp}{\ensuremath{\varnothing}}

\newcommand{\Id}{\ensuremath{\mathrm{Id}}}

\newcommand{\RP}{\ensuremath{\left[0,{+}\infty\right[}}

\newcommand{\NN}{\ensuremath{\mathbb N}}

\newcommand{\weakly}{\ensuremath{\:\rightharpoonup\:}}

\newcommand{\zer}{\ensuremath{\text{\rm zer}\,}}

\newcommand{\proj}{\ensuremath{\text{\rm proj}}}

\newcommand{\gra}{\ensuremath{\text{\rm graph}\,}}

\def\abstract{\noindent{\bfseries Abstract}. \ignorespaces}

\newtheorem{theorem}{Theorem}

\theoremstyle{plain}{\theorembodyfont{\rmfamily}%
}

\theoremstyle{plain}{\theorembodyfont{\rmfamily}%
\newtheorem{cexample}[theorem]{Counterexample}}
\theoremstyle{plain}{\theorembodyfont{\rmfamily}%
}
\theoremstyle{plain}{\theorembodyfont{\rmfamily}%
}
\theoremstyle{plain}{\theorembodyfont{\rmfamily}%
}
\theoremstyle{plain}{\theorembodyfont{\rmfamily}%
}

\let\to\rightarrow

\begin{document}

\title{\sffamily\huge%
The Douglas--Rachford Algorithm Converges Only Weakly\thanks{%
Contact author:
P. L. Combettes.
Email: \email{plc@math.ncsu.edu}.
Phone: +1 919 515 2671.
This work was supported by the National Science
Foundation under grant CCF-1715671.
}
}

\author{Minh N. B\`ui and Patrick L. Combettes\\
\small North Carolina State University,
Department of Mathematics,
Raleigh, NC 27695-8205, USA\\
\small \email{mnbui@ncsu.edu}\: and \:\email{plc@math.ncsu.edu}
}
\date{~}
\maketitle

\begin{abstract}
We show that the weak convergence of the Douglas--Rachford algorithm 
for finding a zero of the sum of two maximally monotone operators 
cannot be improved to strong convergence. Likewise, we show that
strong convergence can fail for the method of partial inverses.
\end{abstract}

\medskip

\begin{keywords}
Douglas--Rachford algorithm, method of partial inverses,
monotone operator, operator splitting, strong convergence.
\end{keywords}

\bigskip

The original Douglas--Rachford splitting algorithm was designed to
decompose positive systems of linear equations \cite{Doug56}.  It
evolved in \cite{Lion79} into a powerful method for finding a zero
of the sum of two maximally monotone operators in Hilbert spaces, a
problem which is ubiquitous in applied mathematics (see
\cite{Livre1} for background on monotone operators). In this
context, the Douglas--Rachford algorithm constitutes a prime
decomposition method in areas such as control, partial differential
equations, optimization, statistics, variational inequalities,
mechanics, optimal transportation, machine learning, and signal
processing. Its asymptotic behavior is described next.

\begin{theorem}
\label{t:dr}
Let $\HH$ be a real Hilbert space, and let $A$ and $B$ be set-valued 
maximally monotone operators from $\HH$ to $2^\HH$ with resolvents 
$J_A=(\Id+A)^{-1}$ and $J_B=(\Id+B)^{-1}$. Suppose that
$\zer(A+B)=\menge{x\in\HH}{0\in Ax+Bx}\neq\emp$, 
let $y_0\in\HH$, and iterate
\begin{equation}
\label{e:dr}
(\forall n\in\NN)\quad
x_n=J_{B}y_n\quad\text{and}\quad
y_{n+1}=y_n+J_{A}(2x_n-y_n)-x_n.
\end{equation}
Then the following hold for some $(y,x)\in\gra J_B$:
\begin{enumerate}
\item
\label{t:dri}
$x=J_A(2x-y)$, $y_n\weakly y$, and $x\in\zer(A+B)$.
\item
\label{t:drii}
$x_n\weakly x$.
\end{enumerate}
\end{theorem}

Property~\ref{t:dri} was established in \cite{Lion79}. Let us note
that, since $J_B$ is not weakly sequentially continuous in general,
the weak convergence of $(y_n)_{n\in\NN}$ in \ref{t:dri} does not 
imply \ref{t:drii}. The latter 
was first established in \cite{Svai11} (see also
\cite[Theorem~26.11(iii)]{Livre1} for an alternate proof). 
In the literature, while various additional conditions on $A$ and $B$
have been proposed to ensure the strong convergence of the sequence
$(x_n)_{n\in\NN}$ in \eqref{e:dr} \cite{Livre1,Joca09,Lion79}, 
it remains an open question whether strong convergence of
$(x_n)_{n\in\NN}$ can fail in the general setting of
Theorem~\ref{t:dr}. We show that this is indeed the case. Our
argument relies on a result of Hundal \cite{Hund04} concerning the 
method of alternating projections.

\begin{cexample}
\label{ex:dr}
In Theorem~\ref{t:dr}, suppose that $\HH$ is infinite-dimensional 
and separable. Let $(e_k)_{k\in\NN}$ be an orthonormal basis of
$\HH$, let $V=\{e_0\}^\bot$, let $y_0=e_2$, and let $K$ be the
smallest closed convex cone containing the set
\begin{equation}
\label{e:K}
\Menge{
e^{-100\xi^3}e_{0}+
\cos\!\big(\pi(\xi-\lfloor\xi\rfloor)/2\big)
e_{\lfloor\xi\rfloor+1}+
\sin\!\big(\pi(\xi-\lfloor\xi\rfloor)/2\big)
e_{\lfloor\xi\rfloor+2}}{\xi\in\RP}.
\end{equation}
Let $\proj_V$ and $\proj_K$ be the projection operators onto $V$ and
$K$, and set 
\begin{equation}
\label{e:AB}
A\colon x\mapsto
\begin{cases}
V^\bot,&\text{if}\;\;x\in V;\\
\emp,&\text{if}\;\;x\notin V
\end{cases}
\qquad\text{and}\quad
B=\big(\proj_V\circ\proj_K\circ\proj_V\big)^{-1}-\Id.
\end{equation}
Then $A$ and $B$ are maximally monotone, and the sequence 
$(x_n)_{\in\NN}$ constructed in Theorem~\ref{t:dr} converges
weakly but not strongly to a zero of $A+B$.
\end{cexample}
\begin{proof}
We first note that $A$ is maximally monotone by virtue of
\cite[Examples~6.43 and and 20.26]{Livre1}. Now set
$T=\proj_V\circ\proj_K\circ\proj_V$. Then it follows from
\cite[Example~4.14]{Livre1} that $T$ is firmly nonexpansive,
that is,
\begin{equation}
(\forall x\in\HH)(\forall y\in\HH)\quad\scal{x-y}{Tx-Ty}\geq
\|Tx-Ty\|^2.
\end{equation}
In turn, we derive from \cite[Proposition~23.10]{Livre1} that
$B=T^{-1}-\Id$ is maximally monotone.
Next, we observe that $0\in\zer A$ and that,
since $K$ is a closed cone,
$0\in K$. Thus, $0=(\proj_V\circ\proj_K\circ\proj_V)0$,
which implies that $0\in\zer B$. Hence,
\begin{equation}
\label{e:0}
0\in \zer(A+B).
\end{equation}
Now set
\begin{equation}
\label{e:z}
z_0=e_2
\quad\text{and}\quad
(\forall n\in\NN)\quad z_{n+1}=\proj_K\big(\proj_Vz_n\big).
\end{equation}
Then, by nonexpansiveness of $\proj_K$,
\begin{equation}
\label{e:8}
(\forall n\in\NN)\quad 
\|z_{n+1}\|^2=\|\proj_K\big(\proj_Vz_n\big)-\proj_K0\|^2
\leq\|\proj_Vz_n\|^2=\|z_n\|^2-\|\proj_Vz_n-z_n\|^2
\end{equation}
and, therefore,
\begin{equation}
\label{e:9}
\proj_Vz_n-z_n\to 0.
\end{equation}
As shown in \cite{Hund04}, we also have
\begin{equation}
\label{e:7}
z_n\weakly 0\quad\text{and}\quad z_n\not\to 0.
\end{equation}
On the other hand, we derive from \eqref{e:AB} that
\begin{equation}
\label{e:J}
J_A=\proj_V\quad\text{and}\quad
J_B=\proj_V\circ\proj_K\circ\proj_V.
\end{equation}
In turn, it follows from \eqref{e:dr} and \eqref{e:z} that
$x_0=\proj_V(\proj_K(\proj_Vz_0))=\proj_Vz_1$.
In addition, $y_0=z_0=\proj_Vz_0$.
Now, assume that, for some $n\in\NN$, $y_n=\proj_Vz_n$
and $x_n=\proj_Vz_{n+1}$. Since $x_n$
and $y_n$ lie in $V$, we derive from
\eqref{e:dr} and \eqref{e:J} that
\begin{equation}
y_{n+1}=y_n+\proj_V(2x_n-y_n)-x_n=x_n=\proj_Vz_{n+1}
\end{equation}
and, in turn, that
\begin{equation}
x_{n+1}=\big(\proj_V\circ\proj_K\circ\proj_V\big)\big(
\proj_Vz_{n+1}\big)=\proj_V\big(\proj_K
\big(\proj_Vz_{n+1}\big)\big)=\proj_Vz_{n+2}.
\end{equation}
Hence,
\begin{equation}
(\forall n\in\NN)\quad x_n=\proj_Vz_{n+1}. 
\end{equation}
Thus, in view of \eqref{e:9}, $x_n-z_{n+1}\to 0$
and we therefore derive from 
\eqref{e:7} and \eqref{e:0} that $x_n\weakly 0\in\zer(A+B)$ and 
$x_n\not\to 0$.
\end{proof}

Next, we settle a similar open question for Spingarn's method of
partial inverses \cite{Spin83} by showing that its strong
convergence can fail.

\begin{theorem}[\rm\cite{Spin83}]
\label{t:spingarn}
Let $\HH$ be a real Hilbert space, let $B\colon\HH\to 2^{\HH}$ be 
maximally monotone, and let $V$ be a closed vector subspace of $\HH$.
Suppose that the problem
\begin{equation}
\label{e:17}
\text{find}\;\;x\in V\;\;\text{and}\;\;u\in V^\bot\;\;
\text{such that}\;\;u\in Bx
\end{equation}
has at least one solution. Let $x_0\in V$,
let $u_0\in V^\bot$, and iterate
\begin{equation}
\label{e:spingarn}
(\forall n\in\NN)\quad
x_{n+1}=\proj_V\big(J_B(x_n+u_n)\big)\quad\text{and}\quad
u_{n+1}=\proj_{V^\bot}\big(J_{B^{-1}}(x_n+u_n)\big).
\end{equation}
Then $(x_n,u_n)_{n\in\NN}$ converges weakly to a solution to
\eqref{e:17}.
\end{theorem}

\begin{cexample}
\label{ex:spingarn}
Define $\HH$, $V$, $K$, and $B$ as in Counterexample~\ref{ex:dr}, 
and set $x_0=e_2$ and $u_0=0$. Then $(0,0)$ solves \eqref{e:17} and 
the sequence $(x_n,u_n)_{n\in\NN}$ constructed in 
Theorem~\ref{t:spingarn} converges weakly but not strongly to 
$(0,0)$.
\end{cexample}
\begin{proof}
Since $J_B=\proj_V\circ\proj_K\circ\proj_V$
and $J_{B^{-1}}=\Id-J_B$, \eqref{e:spingarn} implies that
\begin{equation}
\label{e:22}
(\forall n\in\NN)\quad
\begin{cases}
x_{n+1}=
\big(\proj_V\circ\proj_K\circ\proj_V\big)(x_n+u_n)\\
u_{n+1}=\proj_{V^\bot}\big(x_n+u_n
-\big(\proj_V\circ\proj_K\circ\proj_V\big)(x_n+u_n)\big).
\end{cases}
\end{equation}
It follows that
\begin{equation}
(\forall n\in\NN)\quad
x_{n+1}=\proj_V\big(\proj_Kx_n\big)\quad\text{and}\quad u_n=0.
\end{equation}
Now define $(z_n)_{n\in\NN}$ as in \eqref{e:z}. Then
$(\forall n\in\NN)$ $x_n=\proj_V z_{n}$. Hence, in view of 
\eqref{e:9} and \eqref{e:7}, we conclude that $0\not\!\leftarrow
x_n\!\!\weakly 0$.
\end{proof}

\end{document}